\documentclass[11pt, letter]{article}
\usepackage{beton}
\usepackage{euler}
\usepackage[T1]{fontenc}
\usepackage{amsmath}
\usepackage{amsfonts}
\usepackage[margin=1in]{geometry}
\usepackage{commath}
\usepackage{mathtools}
\usepackage{indentfirst}
\usepackage{amsfonts}
\usepackage{graphicx}
\usepackage{epstopdf}
\usepackage{algorithmic}
\usepackage{xfrac}
\usepackage{subcaption}
\usepackage{hyperref}
\ifpdf
  \DeclareGraphicsExtensions{.eps,.pdf,.png,.jpg}
\else
  \DeclareGraphicsExtensions{.eps}
\fi

\newcommand{\Ce}{\ensuremath{\mathrm{Ce}}}
\newcommand{\ce}{\ensuremath{\mathrm{ce}}}
\newcommand{\Se}{\ensuremath{\mathrm{Se}}}
\newcommand{\se}{\ensuremath{\mathrm{se}}}
\renewcommand{\emph}[1]{\textsl{#1}}
\title{Pure Tone Modes\\
for a $5:3$ Elliptic Drum}
\author{Robert M.~Corless\\
Scientific Director, Ontario Research Centre for Computer Algebra\\
member, The Rotman Institute of Philosophy\\
The University of Western Ontario}
\date{\today}

\begin{document}
\maketitle

\section{Expository Lump}
In 1868 \'Emile Mathieu published a memoir showing how to find the modes, often called \emph{standing modes}, that generate pure tones in the vibration of an elliptic drum of arbitrary aspect ratio~\cite{mathieu1868memoire}. In preparing~\cite{brimacombe2020computation}, together with my co-authors Chris Brimacombe and Mair Zamir, I read that memoir in detail; indeed we asked my colleague Robert H.~C.~Moir to translate the memoir from French in order to make that job easier.  We learned that Mathieu had solved the linear model PDE problem completely, and had given a detailed discussion of the possible nodal lines---that is, places on the membrane that would not move if the drum was oscillating in that mode, including a complete zero-counting discussion using Sturm theory, which he extended to the periodic case.  But he published no pictures; perhaps he had sketched them for himself, but we know of none that survive from that time.  Later authors, of course, have supplied such pictures: see for instance~\cite{gutierrez2003mathieu} or~\cite{chaos2002mathieu} or~\cite{chen1994visualization}.
Once one has seen such pictures, Mathieu's verbal descriptions of the number and location of nodal lines (hyperbolas and ellipses, actually) make more sense.  In this note I show some of the pictures I computed for myself using the (somewhat quixotic) software that I wrote for the Mathieu functions and for solving the Mathieu equation.

We parameterize the elliptic geometry as Mathieu did (he credited the change of coordinates to Lam\'e who called the result ``thermometric'' coordinates) through what is now called the \emph{confocal elliptic coordinates} $x = c\cosh\beta\cos\alpha$, $y = c\sinh\beta\sin\alpha$, or equivalently $x+iy = c\cos(\alpha-i\beta)$.  The real parameter $c$ gives the distance $2c$ between the foci of the ellipse. These coordinates are \emph{singular} in the case $c=0$ when the  ellipse becomes a circle, which is a bit of a headache: the proper analysis of the reduction to the circular case requires a double limit process (see Appendix~B of~\cite{brimacombe2020computation} for a sketch of such an analysis).  Mathieu starts with a wave equation, containing a material parameter $m^2$ representing the ratio of stress to density in the membrane, and then by separation of variables $w = P(\alpha)Q(\beta)\sin 2\lambda m t$ with the frequency of oscillation parameter $\lambda m$ dealing with the physical parameter; here $P$ is a solution of the so-called Mathieu differential equation and $Q$ is a solution of the so-called modified Mathieu equation.  These are also called the \emph{angular} and \emph{radial} Mathieu equations.

The Mathieu differential equation is
\begin{equation}\label{eq:MathieuDE}
  y''(\alpha) + (a - 2q\cos 2\alpha)y(\alpha) = 0
\end{equation}
and the modified Mathieu differential equation is
\begin{equation}\label{eq:ModifiedMathieuDE}
  y''(\beta) - (a - 2q\cosh 2\beta)y(\beta) = 0
\end{equation}
where~$q$ is a parameter that depends on the focal parameter $c$ and the frequency of oscillation $2\lambda m$; the parameter~$a$ is an \emph{eigenvalue} of the differential equation~\eqref{eq:MathieuDE}. The eigenvalue~$a$ depends on the value of~$q = \lambda^2c^2$, and must be chosen to make the solution $y(\alpha)$ periodic with period either~$\pi$ or~$2\pi$.  Mathieu also divided solutions up into \emph{even} and \emph{odd} categories, analogous to cosine and to sine. As in~\cite{brimacombe2020computation} we follow Mathieu and denote the Mathieu functions by $\ce_g(q;\alpha)$ for the even solutions and $\se_g(q;\alpha)$ for the odd, and the modified Mathieu functions by $\Ce_g(q;\beta)$ and $\Se_g(q;\beta)$. Here $g$ is an integer (Mathieu used the letter $g$ in this way, and it might seem unusual to modern eyes that are now used to the \textsc{Fortran} I-N convention of single-letter integer variables being $i$, $j$, $k$, $\ell$, $m$, or $n$). We have $g\ge 0$ for the even solutions and $g\ge 1$ for the odd.  We normalize the eigenfunctions as follows:
\begin{align}\label{eq:normalization}
  \ce_g(q;0) &=  1 \nonumber\\
  \left.\frac{d}{d\alpha}\ce_g(q;\alpha)\right\vert_{\alpha=0} &= 0 \nonumber\\
  \se_g(q;0) &=  0 \nonumber\\
  \left.\frac{d}{d\alpha}\se_g(q;\alpha)\right\vert_{\alpha=0} &= 1\>.
\end{align}
We normalize the modified Mathieu functions analogously, so that $\Ce_g(q;\beta) = \ce_g(q;i\beta)$ and
$\Se_g(q;\beta) = -i\se_g(q;i\beta)$.
This normalization is used in some places in the literature, but most people use a different one involving the square integral over the period.  See~\cite{brimacombe2020computation} for a discussion of why we chose this one.

The frequency of vibration associated with a particular standing mode depends also on the material properties of the membrane as encapsulated in the parameter $m^2$, the ratio of tension to density in the (homogeneous) membrane. If the membrane is oscillating like $\sin 2\lambda m t$ then via $q = \lambda^2 c^2$ we have
\begin{equation}\label{eq:frequency}
  \lambda = \frac{\sqrt{q}}{c} = \frac{h}{c}\>.
\end{equation}
Some authors use $h^2$ in place of~$q$.

It is the so-called \emph{modified} Mathieu functions $\Ce_g(q;\beta) = \ce_g(q;i\beta)$ and $\Se_g(q;\beta) = -i\se_g(q;i\beta)$ which determine, through the boundary conditions (the membrane is fixed at the edges, so $y=0$ there), the value of~$q$.  Mathieu established that every natural mode of the membrane is of the form
\begin{equation}\label{eq:evenmode}
  \Ce_g(q;\beta)\ce_g(q;\alpha)
\end{equation}
or
\begin{equation}\label{eq:oddmode}
  \Se_g(q;\beta)\se_g(q;\alpha)\>.
\end{equation}
The dependence on the frequency of oscillation can be made clearer by writing a mode in full as
\begin{equation}\label{eq:witht}
  P(q;\alpha)Q(q;\beta)\sin2\lambda m t = P(\lambda^2c^2;\alpha)Q(\lambda^2c^2;\beta)\sin2\lambda m t\>.
\end{equation}
That is, \emph{multiple values of}~$q$ are possible for the same PDE with fixed parameter $m$ because different modes of vibration can be excited at the same time.  Therefore a general solution composed of a sum of different excited modes would be a sum over modes with different values of~$\lambda$ and therefore of~$q$.

If we have a fixed elliptical drum in mind, then it will have a given aspect ratio or equivalently eccentricity.  Say, for the sake of argument, that the aspect ratio is $5:3$.  Then we have (putting the major axis on the $x$-axis)
\begin{align*}
  5 &= c\cosh\beta_0 \\
  3 &= c\sinh\beta_0
\end{align*}
defining the outer boundary of the ellipse (this also defines the units we use): this is quickly solved\footnote{$\tanh\beta_0 = 3/5$ so $\beta_0 = \ln 2$ and $5^2-3^2 = c^2 = 16$.} to get $c=4$ and $\beta_0=\ln 2$.  Part of the attraction of choosing this aspect ratio for demonstration is that the numbers are so simple.

\begin{figure}
  \centering
  \subcaptionbox{$\Ce_0(1.7353;\beta)=\ce_0(1.7353;i\beta)$\label{fig:modce0}}{\includegraphics[width=6cm]{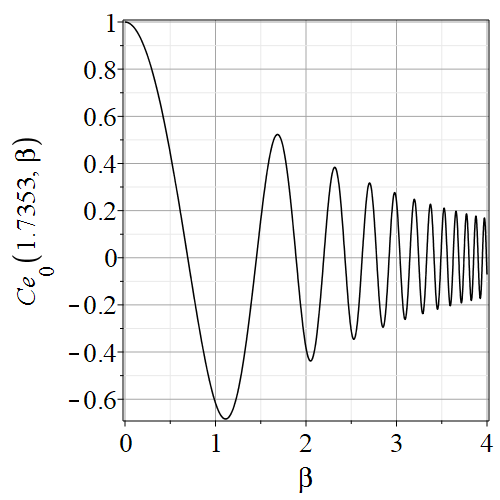}}
  \subcaptionbox{$\Se_1(5.4300;\beta)=\Im(\se_1(5.4300;i\beta))$\label{fig:modse1}}{\includegraphics[width=6cm]{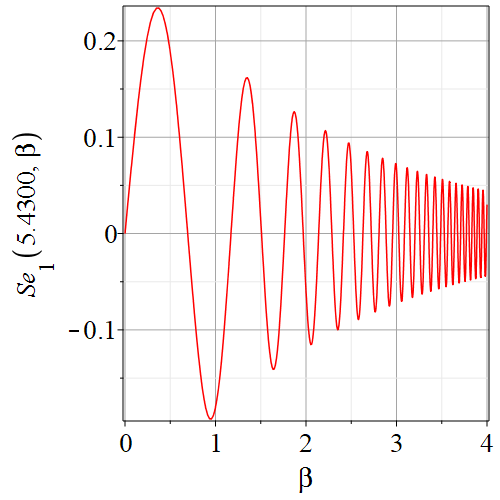}}
  \caption{(Left) A graph of $\Ce_0(q;\beta)$ when $q=1.7353$. As the argument $\beta$ increases, the function becomes increasingly oscillatory.  This value of~$q$ could be used for an elliptical drum whose vertical dimension was such as to coincide with a zero of this function (units depending on the locations of the foci at $\pm c$). By construction, the first such zero occurs at $\beta = 0.6931$, giving an aspect ratio of $5:3$.
  (Right) A graph of $\Se_1(q;\beta)=\Im(\se_1(q;i\beta))$ when $q=5.4300$. Again the first zero allows a pure tone for an ellipse of aspect ratio $5:3$, by construction.}\label{fig:modifiedMathieu0q2}
\end{figure}

To find the permissible values of~$q$ that leave the edge of the drum fixed, then, we have to solve either
\begin{equation}\label{eq:Ceeq}
  \Ce_g(q;\ln2) = 0
\end{equation}
or
\begin{equation}\label{eq:Ceeq}
  \Se_g(q;\ln2) = 0
\end{equation}
for the parameter~$q$.  These are (of course) nonlinear equations, and for real values of~$\lambda$ and $c$ turn out always to have an infinite number of solutions, for each value of the index $g$.  See figure~\ref{fig:modifiedMathieu0q2} where we show a representative pair of modified Mathieu equations and see figure~\ref{fig:qln2} where we fix $\beta=\beta_0 = \ln2$ and vary $q$, exhibiting an oscillatory function which apparently has an infinite number of zeros.  In practice we are only interested in the first few; the very large values of~$q$ giving distant zeros will give impractically high-frequency standing modes anyway.

\section{Finding~$q$ numerically}
\begin{figure}
  \centering
  \subcaptionbox{$\Ce_1(q;\ln2)$\label{fig:Ce1qln2}}{\includegraphics[width=4cm]{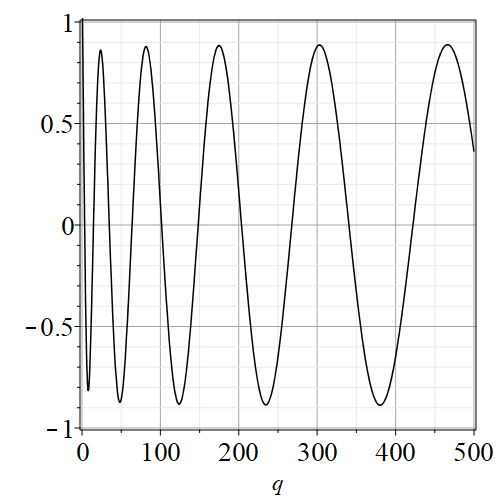}}
  \subcaptionbox{$2\sqrt{q}\Se_1(q;\ln2)$\label{fig:Se1qln2}}{\includegraphics[width=4cm]{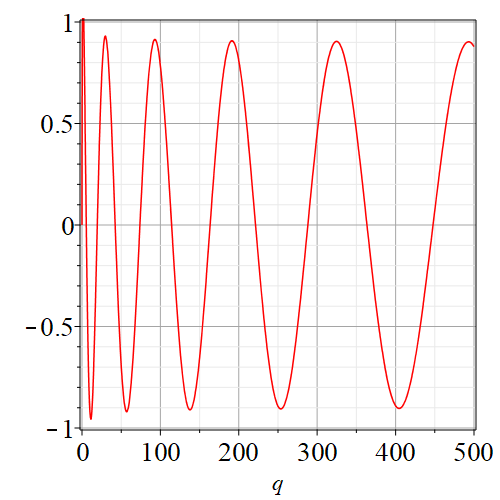}}
  \subcaptionbox{$\Ce_2(q;\ln2)$\label{fig:Ce2qln2}}{\includegraphics[width=4cm]{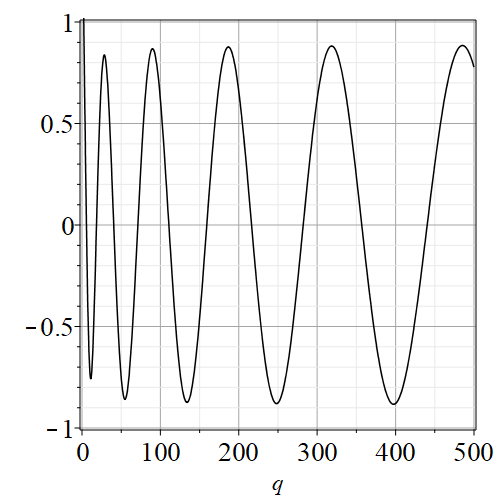}}
  \subcaptionbox{$2\sqrt{q}\Se_2(q;\ln2)$\label{fig:Se2qln2}}{\includegraphics[width=4cm]{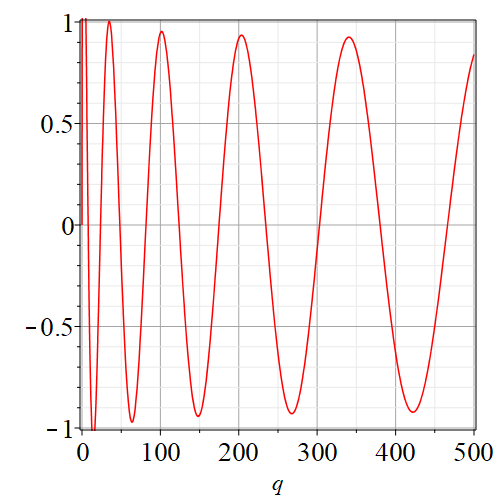}}
  \caption{Modified Mathieu functions for $\beta=\ln2$ as a function of $q$.  From the graphs, we would find it natural that there are an infinite number of values of $q$ for which $\Ce_g(q;\ln2)=0$ and similarly $\Se_g(q;\ln2)=0$. With the normalization used in this paper, it seems that $\Ce_g(q;\ln2)$ eventually oscillates between $-1$ and $1$, while $\Se_g(q;\ln2)$ decays like $1/(2\sqrt{q})$. }\label{fig:qln2}
\end{figure}

We wish to find, for a given nonnegative integer $g$ and choice of ``even'' or ``odd'', the values of~$q$ for which $Q(q;\ln 2) = 0$, where $Q(q;\beta)$ is either $\Ce_g(q;\beta)$ (even case) or $\Se_g(q;\beta)$ (odd case) for the appropriate $g$.  Either way this is a nonlinear equation in the real parameter~$q$ that has an infinite number of solutions, as evidenced by figures~\ref{fig:Ce1qln2}--\ref{fig:Se2qln2}.  In those figures, it is important that we chose to normalize by making $\Ce_g(q;0)=1$ while its derivative with respect to $\beta$ was zero there; if instead we had chosen to use the integral normalization then these values seem to decay root-exponentially, i.e. like $\exp(-2\sqrt{q})$.

We also wish to find~$q$ so that $\ln2$ is the \emph{first} zero of the function $Q(q;\beta)$, and then another value of~$q$ so that $\ln 2$ is the \emph{second} zero, and so on. Each choice of $g$ and of ``even'' versus ``odd' will give a countably infinite number of values of~$q$.  There is no standard notation for these zeros. For now we will write $q_{g,k}^e$ for the for the value of~$q$ that makes $\ln 2$ the $k$th zero of $\Ce_g(q;\beta)$ (the even case), for $k = 1$, $2$, $3$, $\ldots$.  Similarly we will write $q_{g,k}^o$ for the value of~$q$ that makes $\ln 2 $ the $k$th zero of $\Se_g(q;\beta)$ (note that $0 = \Se_g(q;0)$ for all~$q$ and we can say this is the $0$th zero, if we like: then here as with $\Ce$ the index $k$ starts at $1$).

By trying to plot $Q(q;\beta)$ as a function of $\beta$ for various~$q$ (we need first to find the appropriate eigenvalue $a$, and how to do this is described in~\cite{brimacombe2020computation} for example) we see that in general increasing~$q$ brings the first zero closer to the origin.  However, the motion of the zeros is not, in general, monotonic.  Nonetheless, by simply plotting the function for a few values of~$q$ one can quickly get an initial estimate $q_0$ of the value of~$q$ for which $Q(q;\ln2) = 0$.

Once one has a decent initial estimate $q_0$, one would like to use Newton's method to find the value of~$q$ more precisely.  This requires derivatives with respect to~$q$, which ordinarily would require a Fr\'echet derivative and the use of a Green's function (which works perfectly well) but I can report with pleasure that the simple Squire-Trapp formula~\cite{Squire(1998)} for numerical differentiation works beautifully here (see~\cite[p.~479]{corless2013graduate} for a discussion of the excellent numerical properties of this formula).  The derivation of the formula is simple: if $f$ is analytic and $x$ and $h$ are real, then
\begin{equation}\label{eq:Taylor}
  f(x + ih) = f(x) + ih f'(x) - \tfrac12h^2f''(x) -i\,\tfrac16 h^3 f'''(x) + O(h^4)\>.
\end{equation}
Therefore, taking the imaginary part gives
\begin{align}\label{eq:SquireTrapp}
  \frac{\Im(f(z+ih))}{h} &= f'(x) - \tfrac16 h^2 f'''(x) + O(h^4)\>, \nonumber\\
                         &= f'(x) + O(h^2)\>.
\end{align}
So all we need to do to compute $f'(x)$ is to compute $f$ near to $x$ but with a small \emph{purely imaginary} perturbation.
The formula does require $f$ to be analytic because it effectively uses the Cauchy-Riemann equations (here $Q$ is in fact entire) and requires $f$ to be real when $x$ is real and to be computed to high relative accuracy in both the real and imaginary part.  This holds in our case.

That is, instead of computing with our \emph{real} value of~$q$, we choose an $h$ smaller than the square root of machine epsilon (or whatever tolerance we are computing $Q(q;\beta)$ to) and compute instead with $q + ih$; we use the real part of the answer for $Q(q;\ln2)$ and the imaginary part (divided by $h$) as $dQ/dq$.  Since $h^2$ will then be smaller than machine epsilon, we can ignore the effects of the perturbation on the real part.  This is a kind of finite difference formula, but one that for analytic $f$ suffers no instability as $h \to 0$.  Really, it's ridiculously effective.

We may then use Newton's method
\[
q_{n+1} = q_n - \frac{Q(q_n;\ln2)}{Q'(q_n;\ln2)}
\]
where the prime now denotes differentiation with respect to~$q$.  All of the results in tables~\ref{tab:qa} and~\ref{tab:qb} were computed in this way.  I only report the values of~$q$ to a few figures; more can be computed on demand.

The \emph{eigenvalues} $a_g(q_n)$ and $b_g(q_n)$ were computed with an even faster iteration, the Inverse Cubic Iteration described in~\cite{corless2020inverse}, essentially just for fun.  The results are not tabulated here because they can be recomputed easily, given~$q$.
\begin{table}[h!]
  \centering
  \caption{Values of the parameter~$q$ for which $\Ce_g(q,\ln2)=0$.  The $i$th table column gives the values of~$q$ for which the $i$th zero of $\Ce_g(q;\beta)$ is $\beta=\ln2$. The corresponding eigenvalues are $a_g(q)$.}\label{tab:qa}
\begin{tabular}{c|r|r|r|r}
   $g$ & first & second & third & fourth \\
   \hline
   0 & 1.7353& 11.356& 29.795& 57.011
\\  1 & 3.3522& 14.627& 34.844& 63.848
\\  2 & 5.6530& 18.486& 40.457& 71.241
\\  3 & 8.6576& 22.968& 46.658& 79.201
\\  4 &12.368& 28.100& 53.463& 87.754
\\  5 &16.779& 33.913& 60.891& 96.903
\end{tabular}
\end{table}
\begin{table}[h!]
  \centering
  \caption{Values of the parameter~$q$ for which $\Se_g(q,\ln2)=0$.  The $i$th table column gives the values of~$q$ for which the $i$th zero of $\Se_g(q;\beta)$ is $\beta=\ln2$. The corresponding eigenvalues are $b_g(q)$.}\label{tab:qb}
\begin{tabular}{c|r|r|r|r}
   $g$ & first & second & third & fourth \\
   \hline
1 & 5.4300& 19.478& 42.306& 73.902
\\ 2 & 7.8189& 23.636& 48.248& 81.626
\\ 3 & 10.803& 28.363& 54.775& 89.914
\\ 4 &14.406& 33.696& 61.800& 98.762
\\ 5 & 18.637& 39.642& 69.499& 108.17
\end{tabular}
\end{table}
\section{Nodal lines}
Once one has the values of~$q$ which are needed, it is straightforward to compute the mode shapes.  Somewhat surprisingly, it is awkward to plot them (either in Maple or in Matlab) because the surface is defined parametrically; most contour plotting software wants a regular $x$-$y$ grid to work from.  In Maple, this can be worked around by using \texttt{plots[surfdata]} which will quickly do a contour plot of irregular data as a decoration of its 3D plots, which can be looked at from the top down (use the keyword \texttt{orientation=[-90,0]} in the \texttt{plots[surfdata]} command).  Unfortunately, this is somewhat ``too fancy'' and does not produce a nice-looking simple 2D contour plot.

Instead, I inverted the coordinate transformation, so $\alpha - i\beta = \arccos( (x+iy)/4 )$, and simply evaluated my numerical solutions of the Mathieu equation (represented as strings of blends~\cite{corless2020blends}) at the resulting~$\alpha$ and~$\beta$.  At low resolution the contour plots take only a few seconds for each one, but take up to about a minute for high resolutions such as a $400$ by $400$ grid.  In contrast, the \texttt{plots[surfdata]} approach only takes about three seconds in total for a similar high-resolution graph.  But in the meantime, the code was good enough to allow the production of the contour plots in this paper.

\begin{figure}
  \centering
  \subcaptionbox{Contour Plot\label{fig:graycontour}}{\includegraphics[width=6cm]{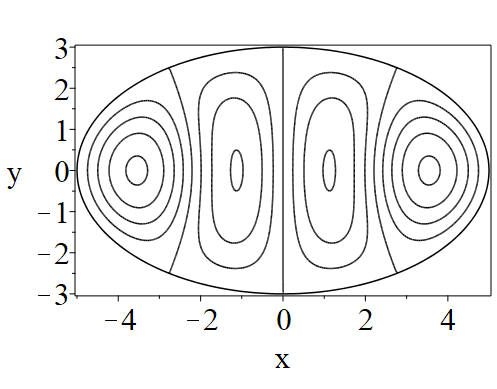}}
  \subcaptionbox{\texttt{plots[surfdata]}\label{fig:colourcontour3}}{\includegraphics[width=6cm]{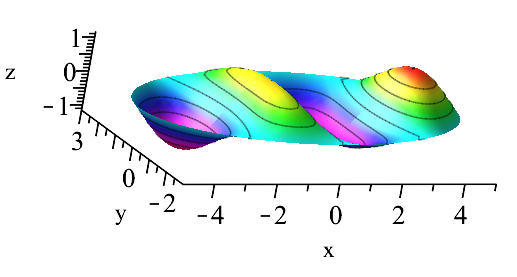}}
  \caption{A typical pure tone mode of a $5:3$ ellipse. Both graphs show $\Ce_3(q;\beta)\ce_3(q;\alpha)$ where $q=8.6576$, $x=4\cosh\beta\cos\alpha$, $y=4\sinh\beta\sin\alpha$. The parameters $-\pi < \alpha < \pi$ and $0 \le \beta \le \ln2$. At $\beta=\ln2$, $4\cosh\beta = 5$ and $4\sinh\beta = 3$, delimiting the elliptic boundary. The left figure is a full version of figure~\ref{fig:a3z1}. The right shows a colour 3D plot of the same.}\label{fig:fullellipse}
\end{figure}

\begin{figure}[t!]
\centering     
\subcaptionbox{$q=1.7353$\label{fig:a0z1} }{\includegraphics[width=40mm]{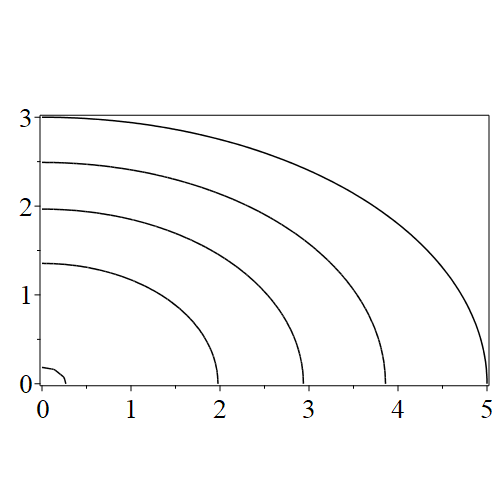}}
\subcaptionbox{$q=11.356$\label{fig:a0z2}}{\includegraphics[width=40mm]{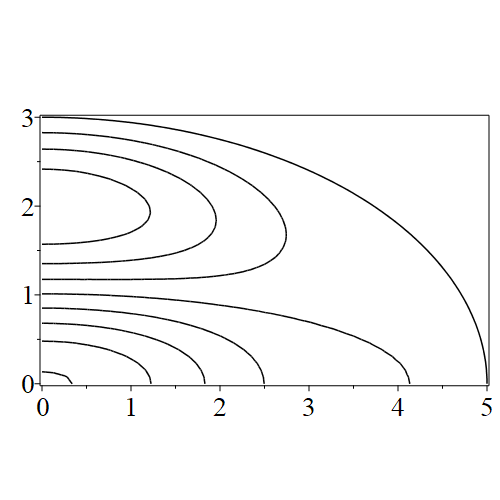}}
\subcaptionbox{$q=29.795$\label{fig:a0z3}}{\includegraphics[width=40mm]{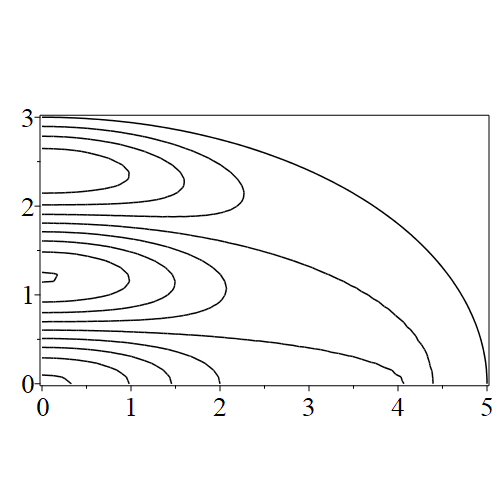}}
\subcaptionbox{$q=57.011$\label{fig:a0z4}}{\includegraphics[width=40mm]{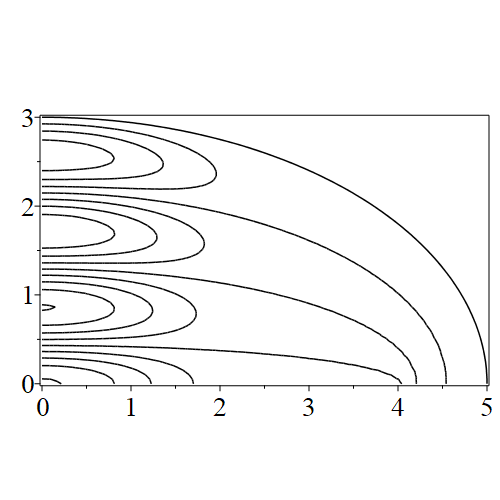}}
\caption{Contours of $\Ce_0(q;\beta)\ce_0(q;\alpha)$ in the first quadrant.
\label{fig:a0}}
\end{figure}

\begin{figure}[t!]
\centering     
\subcaptionbox{$q=3.3522$\label{fig:a1z1} }{\includegraphics[width=40mm]{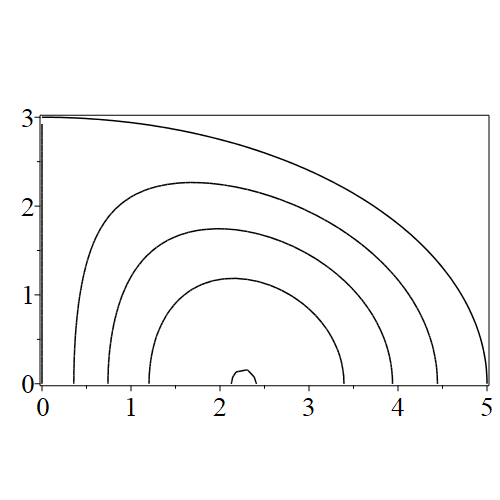}}
\subcaptionbox{$q=14.627$\label{fig:a1z2}}{\includegraphics[width=40mm]{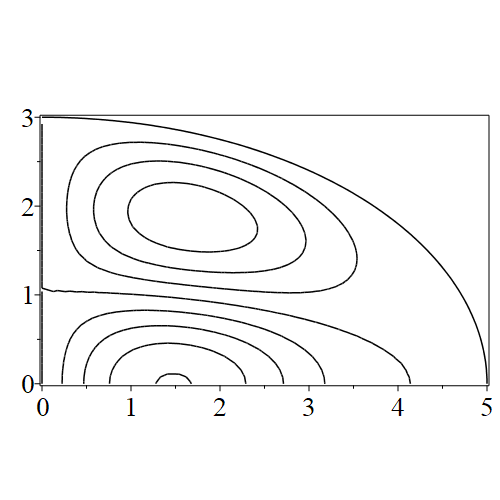}}
\subcaptionbox{$q=34.844$\label{fig:a1z3}}{\includegraphics[width=40mm]{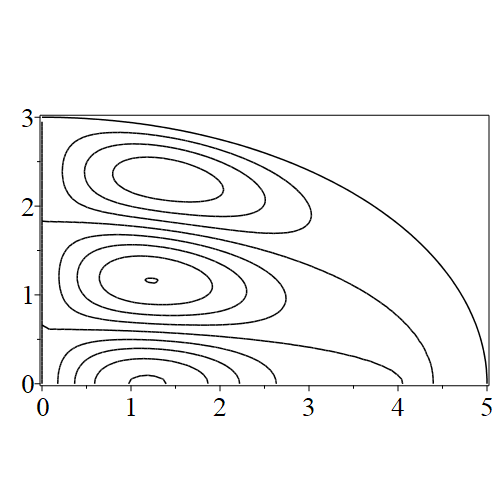}}
\subcaptionbox{$q=63.848$\label{fig:a1z4}}{\includegraphics[width=40mm]{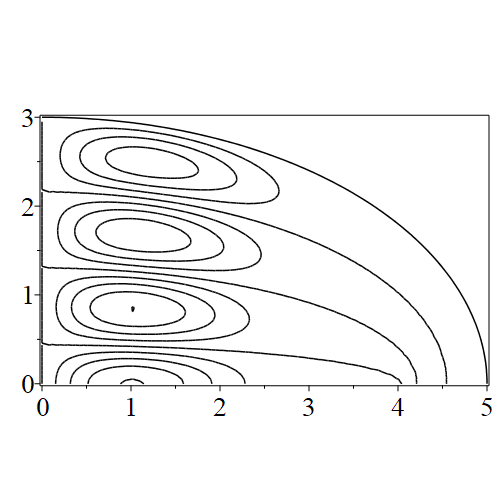}}
\caption{Contours of $\Ce_1(q;\beta)\ce_1(q;\alpha)$ in the first quadrant.
\label{fig:a1}}
\end{figure}

\begin{figure}[t!]
\centering     
\subcaptionbox{$q=5.4300$\label{fig:b1z1} }{\includegraphics[width=40mm]{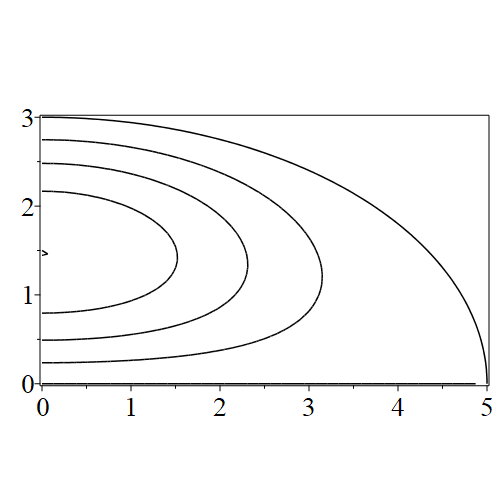}}
\subcaptionbox{$q= 19.478$\label{fig:b1z2}}{\includegraphics[width=40mm]{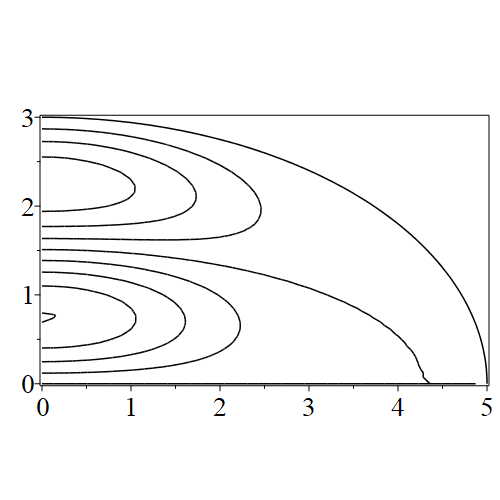}}
\subcaptionbox{$q= 42.306$\label{fig:b1z3}}{\includegraphics[width=40mm]{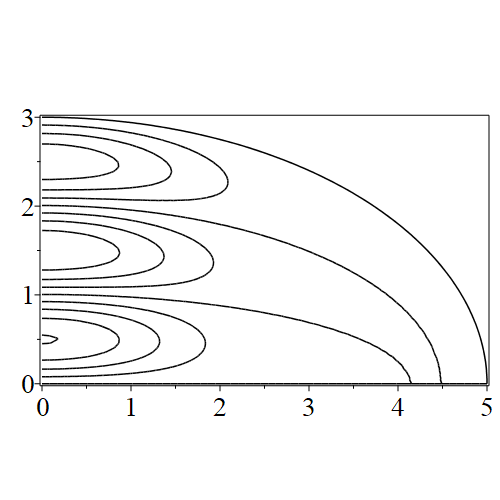}}
\subcaptionbox{$q= 73.902$\label{fig:b1z4}}{\includegraphics[width=40mm]{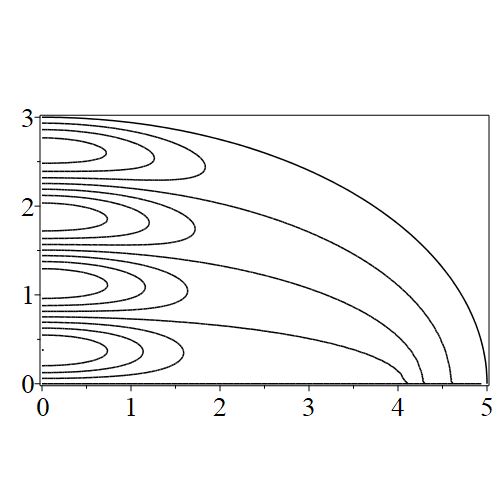}}
\caption{Contours of $\Se_1(q;\beta)\se_1(q;\alpha)$ in the first quadrant.
\label{fig:b1}}
\end{figure}

\begin{figure}[t!]
\centering     
\subcaptionbox{$q=5.6530$\label{fig:a2z1} }{\includegraphics[width=40mm]{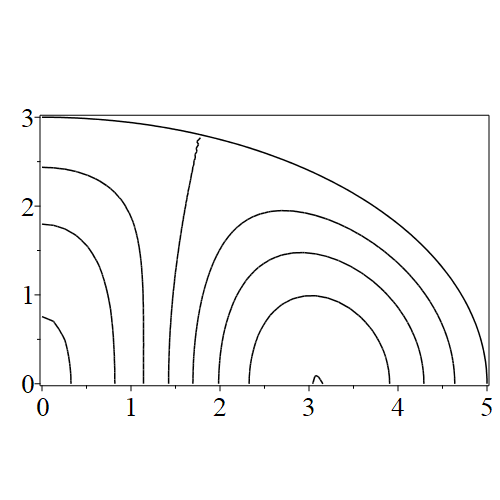}}
\subcaptionbox{$q= 18.486$\label{fig:a2z2}}{\includegraphics[width=40mm]{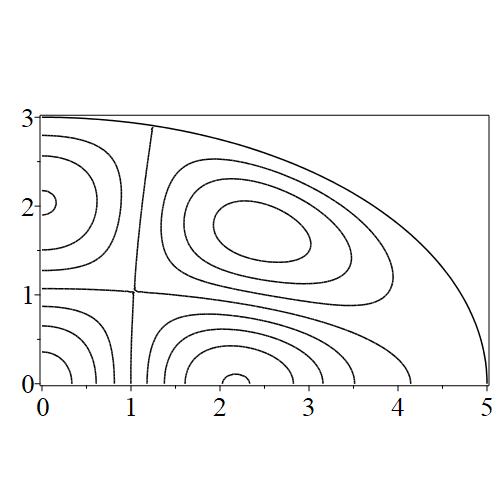}}
\subcaptionbox{$q= 40.457$\label{fig:a2z3}}{\includegraphics[width=40mm]{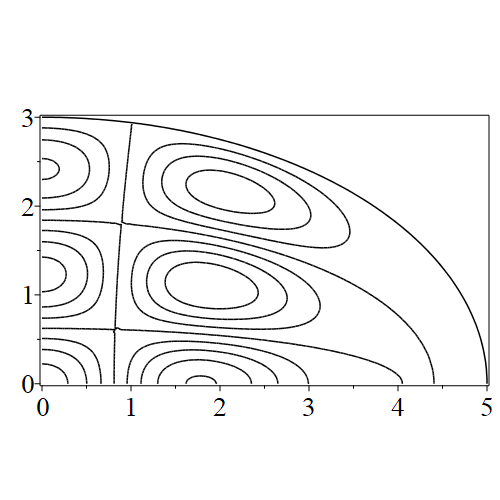}}
\subcaptionbox{$q= 71.241$\label{fig:a2z4}}{\includegraphics[width=40mm]{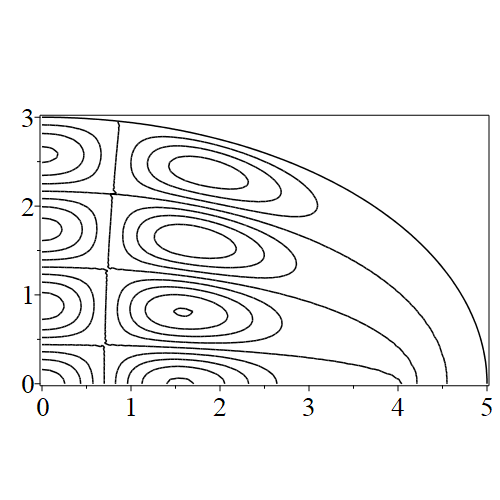}}
\caption{Contours of $\Ce_2(q;\beta)\ce_2(q;\alpha)$ in the first quadrant.
\label{fig:a2}}
\end{figure}

\begin{figure}[t!]
\centering     
\subcaptionbox{$q=7.8189$\label{fig:b2z1} }{\includegraphics[width=40mm]{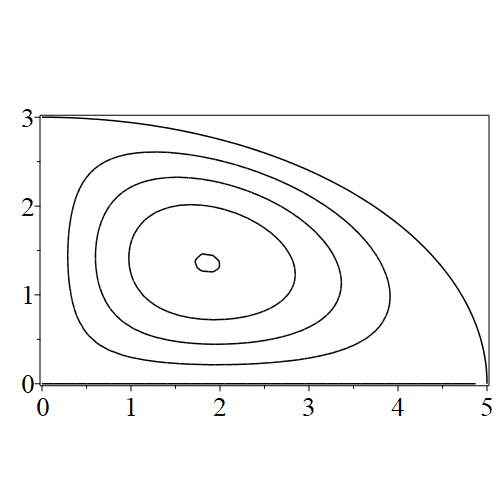}}
\subcaptionbox{$q= 23.636$\label{fig:b2z2}}{\includegraphics[width=40mm]{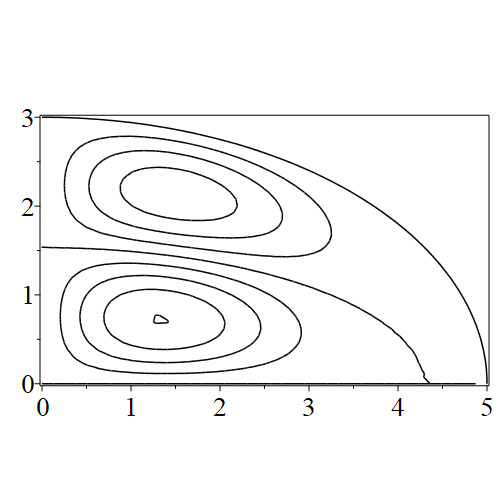}}
\subcaptionbox{$q= 48.248$\label{fig:b2z3}}{\includegraphics[width=40mm]{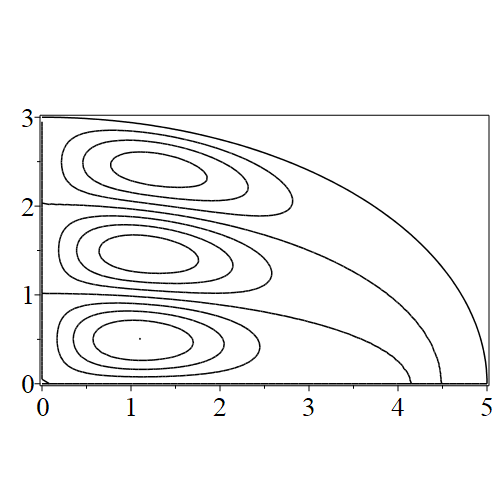}}
\subcaptionbox{$q= 81.626$\label{fig:b2z4}}{\includegraphics[width=40mm]{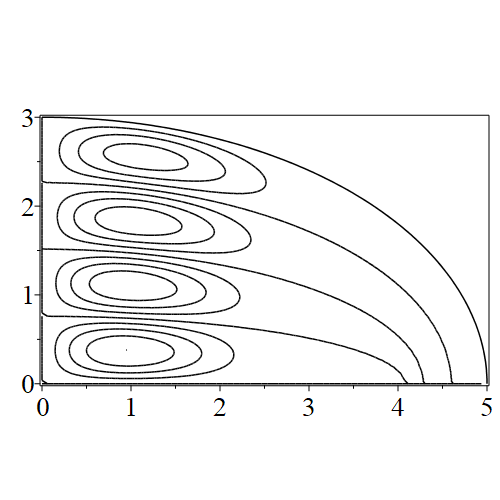}}
\caption{Contours of $\Se_2(q;\beta)\se_2(q;\alpha)$ in the first quadrant.
\label{fig:b2}}
\end{figure}

\begin{figure}[t!]
\centering     
\subcaptionbox{$q=8.6576$\label{fig:a3z1} }{\includegraphics[width=40mm]{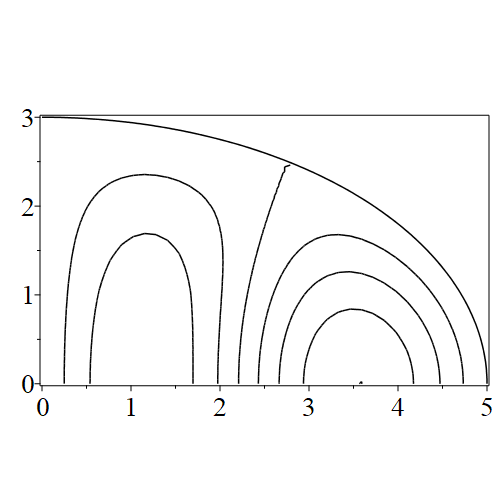}}
\subcaptionbox{$q= 22.968$\label{fig:a3z2}}{\includegraphics[width=40mm]{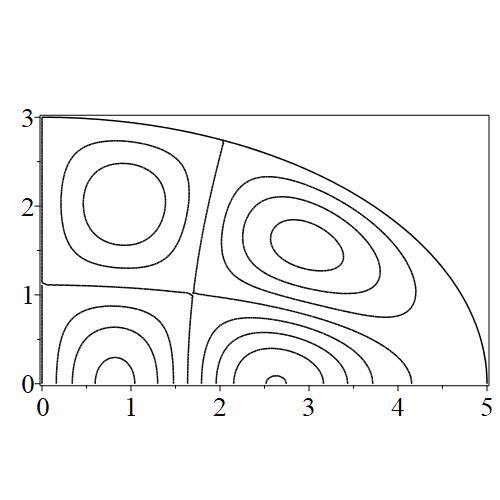}}
\subcaptionbox{$q= 46.658$\label{fig:a3z3}}{\includegraphics[width=40mm]{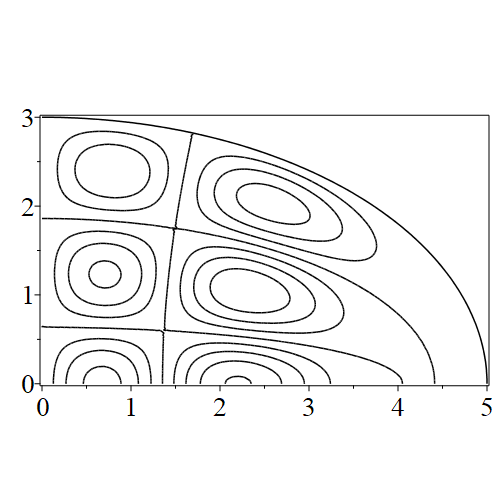}}
\subcaptionbox{$q= 79.201$\label{fig:a3z4}}{\includegraphics[width=40mm]{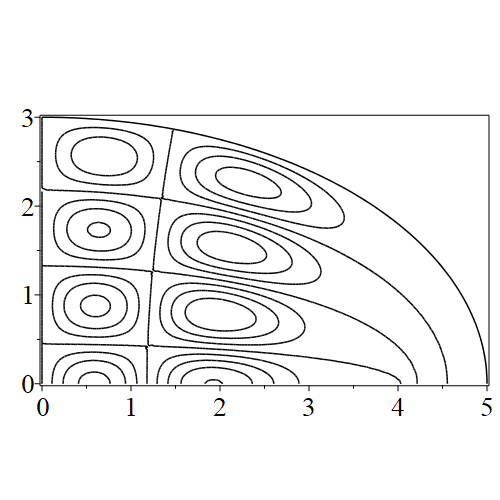}}
\caption{Contours of $\Ce_3(q;\beta)\ce_3(q;\alpha)$ in the first quadrant.
\label{fig:a3}}
\end{figure}

\begin{figure}[t!]
\centering     
\subcaptionbox{$q=10.803$\label{fig:b3z1} }{\includegraphics[width=40mm]{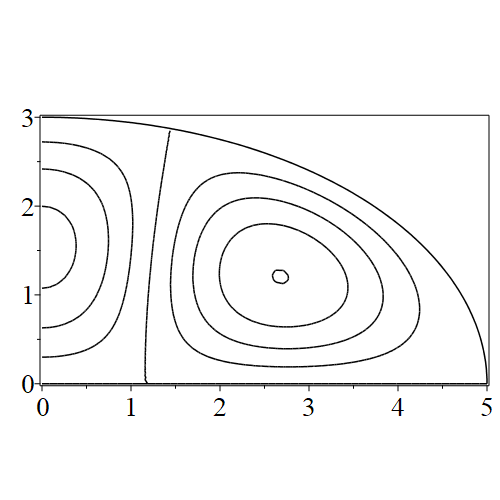}}
\subcaptionbox{$q= 28.363$\label{fig:b3z2}}{\includegraphics[width=40mm]{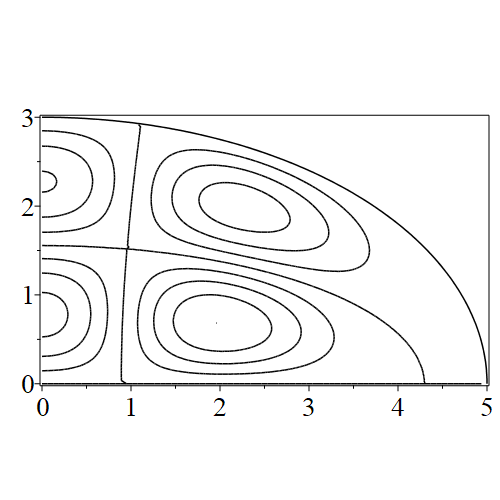}}
\subcaptionbox{$q= 54.775$\label{fig:b3z3}}{\includegraphics[width=40mm]{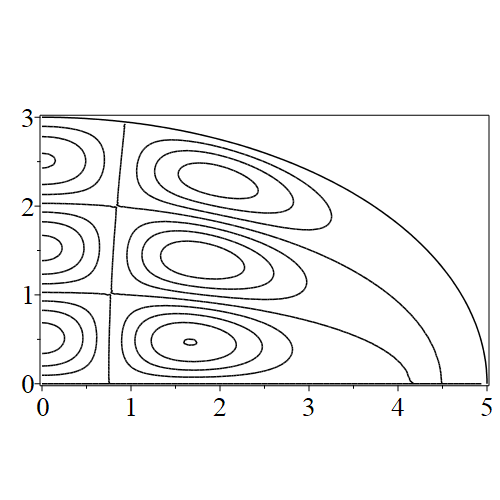}}
\subcaptionbox{$q= 89.914$\label{fig:b3z4}}{\includegraphics[width=40mm]{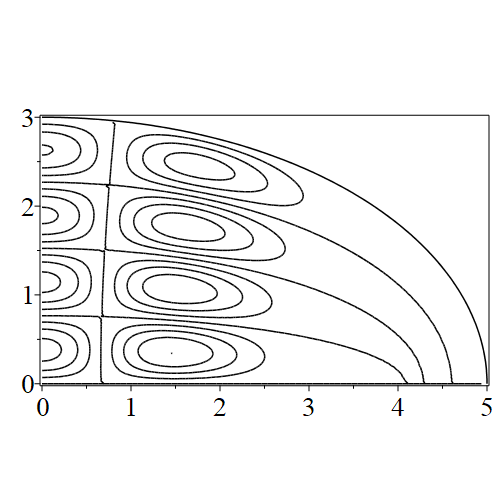}}
\caption{Contours of $\Se_3(q;\beta)\se_3(q;\alpha)$ in the first quadrant.
\label{fig:b3}}
\end{figure}

\begin{figure}[t!]
\centering     
\subcaptionbox{$q=12.368$\label{fig:a4z1} }{\includegraphics[width=40mm]{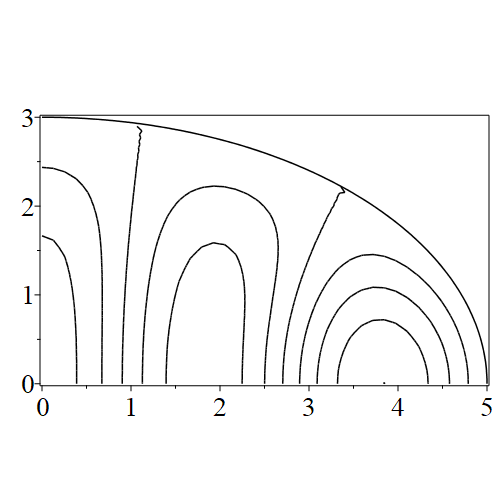}}
\subcaptionbox{$q= 28.100$\label{fig:a4z2}}{\includegraphics[width=40mm]{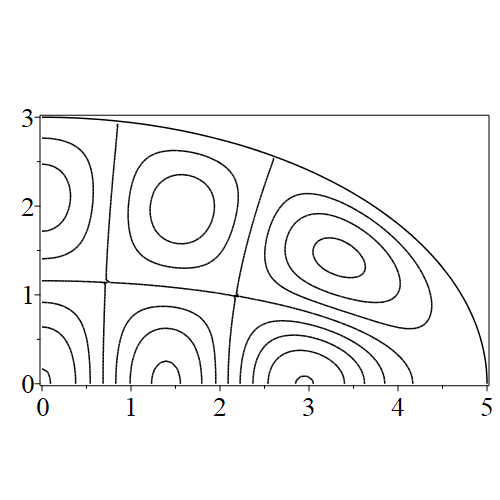}}
\subcaptionbox{$q= 53.463$\label{fig:a4z3}}{\includegraphics[width=40mm]{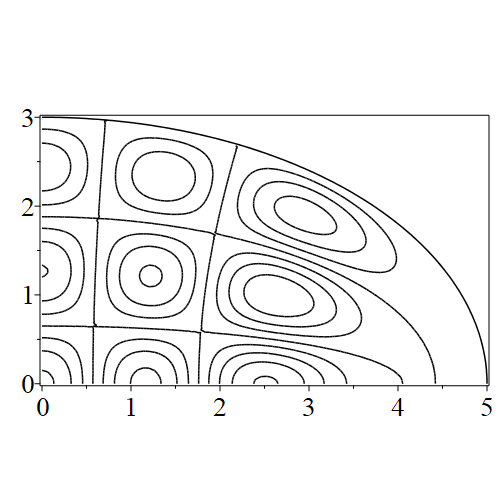}}
\subcaptionbox{$q= 87.754$\label{fig:a4z4}}{\includegraphics[width=40mm]{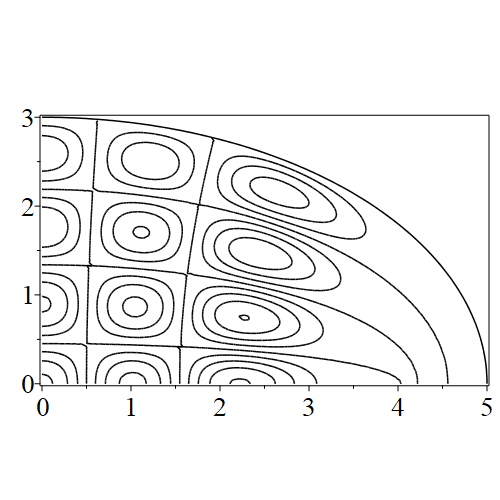}}
\caption{Contours of $\Ce_4(q;\beta)\ce_4(q;\alpha)$ in the first quadrant.
\label{fig:a4}}
\end{figure}

\begin{figure}[t!]
\centering     
\subcaptionbox{$q=14.406$\label{fig:b4z1} }{\includegraphics[width=40mm]{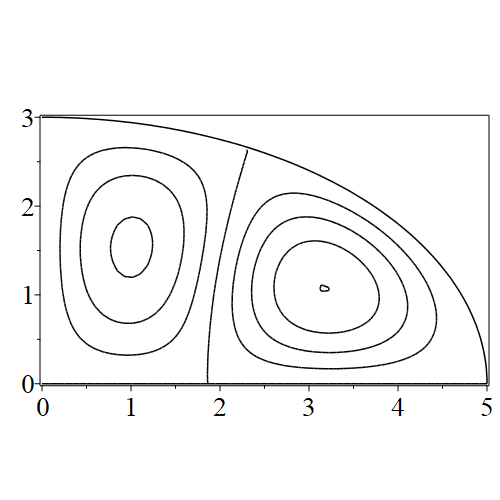}}
\subcaptionbox{$q= 33.696$\label{fig:b4z2}}{\includegraphics[width=40mm]{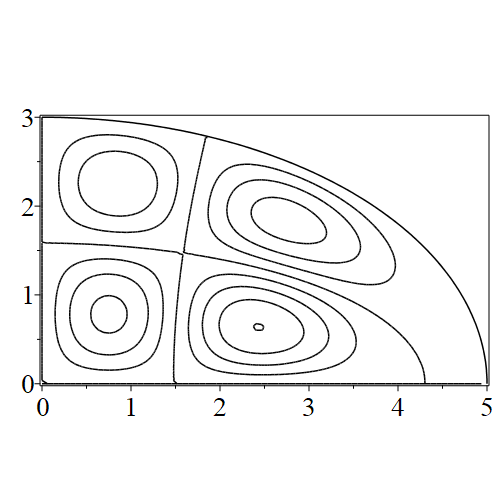}}
\subcaptionbox{$q= 61.800$\label{fig:b4z3}}{\includegraphics[width=40mm]{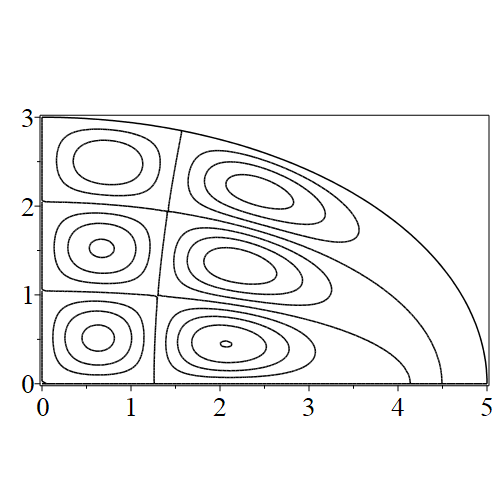}}
\subcaptionbox{$q= 98.762$\label{fig:b4z4}}{\includegraphics[width=40mm]{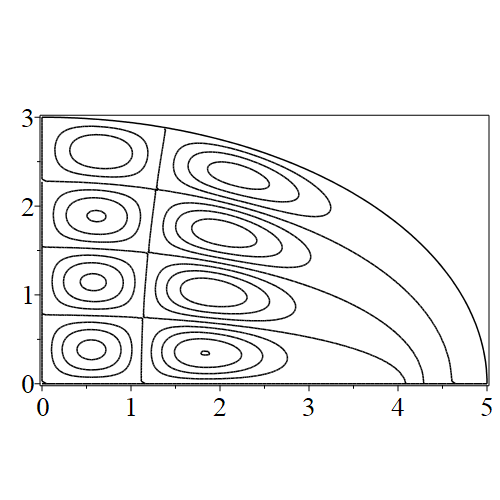}}
\caption{Contours of $\Se_4(q;\beta)\se_4(q;\alpha)$ in the first quadrant.
\label{fig:b4}}
\end{figure}

\begin{figure}[t!]
\centering     
\subcaptionbox{$q=16.779$\label{fig:a5z1} }{\includegraphics[width=40mm]{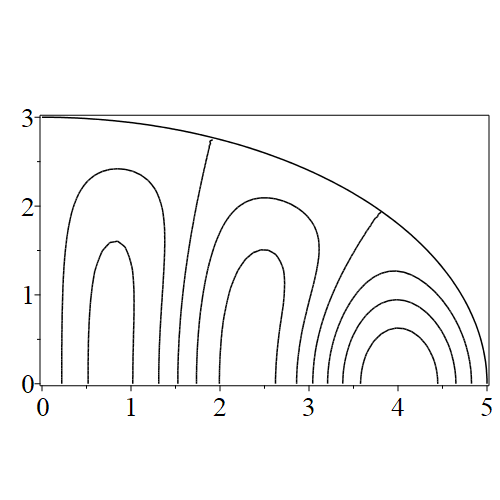}}
\subcaptionbox{$q= 33.913$\label{fig:a5z2}}{\includegraphics[width=40mm]{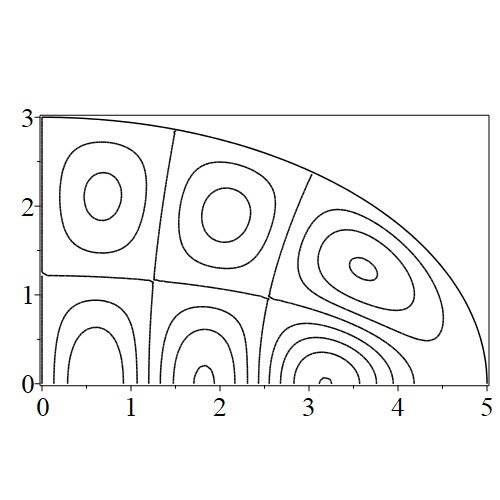}}
\subcaptionbox{$q= 60.891$\label{fig:a5z3}}{\includegraphics[width=40mm]{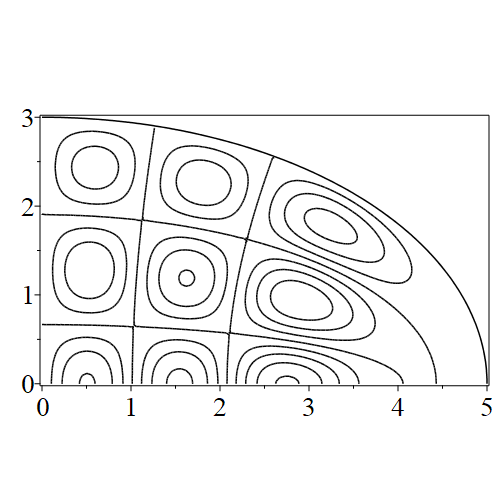}}
\subcaptionbox{$q= 96.903$\label{fig:a5z4}}{\includegraphics[width=40mm]{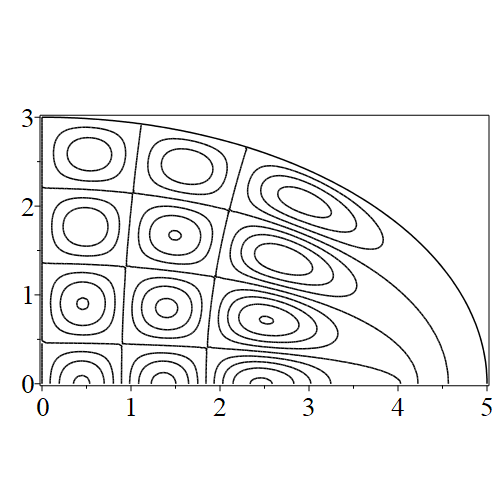}}
\caption{Contours of $\Ce_5(q;\beta)\ce_5(q;\alpha)$ in the first quadrant.
\label{fig:a5}}
\end{figure}

\begin{figure}[t!]
\centering     
\subcaptionbox{$q=18.637$\label{fig:b5z1} }{\includegraphics[width=40mm]{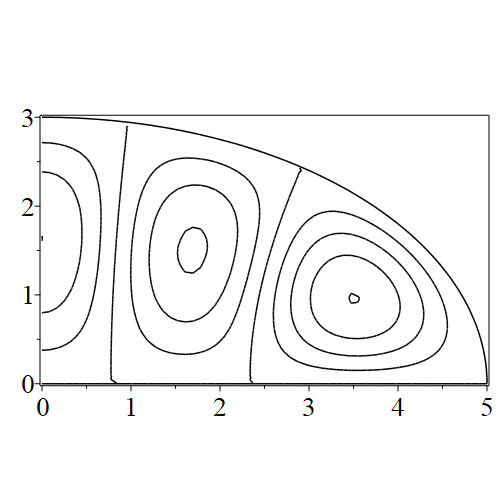}}
\subcaptionbox{$q= 39.642$\label{fig:b5z2}}{\includegraphics[width=40mm]{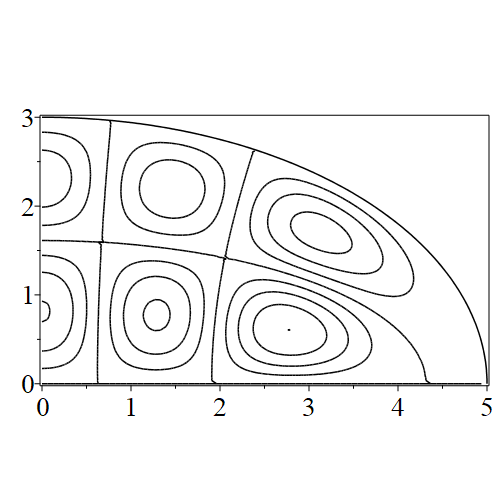}}
\subcaptionbox{$q= 69.499$\label{fig:b5z3}}{\includegraphics[width=40mm]{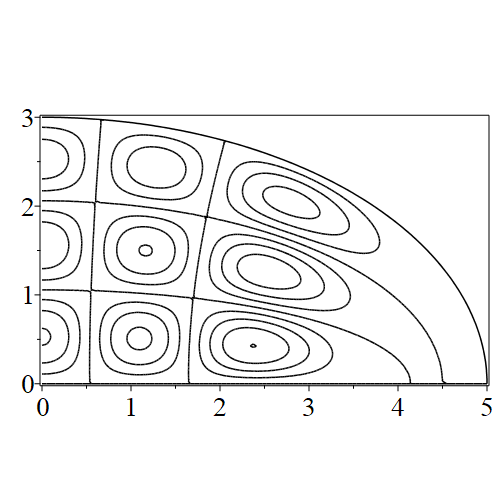}}
\subcaptionbox{$q= 108.17$\label{fig:b5z4}}{\includegraphics[width=40mm]{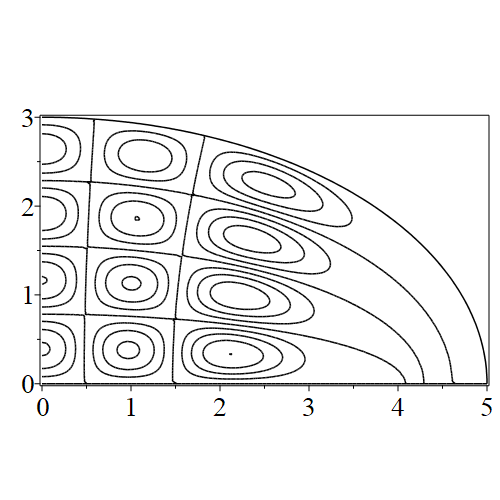}}
\caption{Contours of $\Se_5(q;\beta)\se_5(q;\alpha)$ in the first quadrant.
\label{fig:b5}}
\end{figure}

In figure~\ref{fig:fullellipse} we see the contours corresponding to an even solution with $g=3$, where the value of~$q$ was chosen so that $\ln2$ was the \emph{first} zero of $\Ce_3(q;\beta)$.  Figures~\ref{fig:a0}--\ref{fig:b5} plot only the contours of the first quadrant; by symmetry the other three can be deduced.  We see confocal elliptic nodal lines, and hyperbolic nodal lines, in several figures; by comparing all the figures we may understand what Mathieu was talking about when he discussed the natural modes of vibration of an elliptic drum.

\begin{figure}[b]
  \centering
  \subcaptionbox{$g=11$\label{fig:whisper11}}[6cm][l]{\includegraphics[width=6cm]{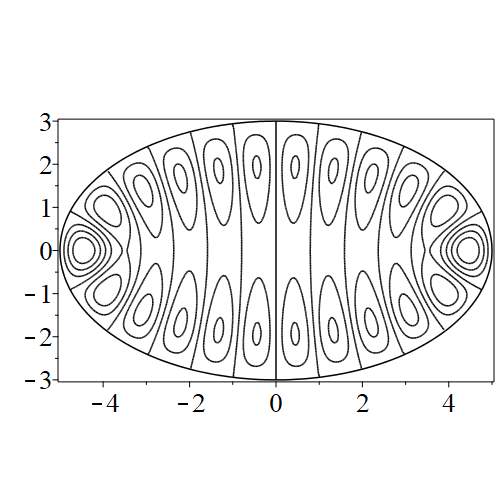}}
  \subcaptionbox{$g=13$\label{fig:whisper13}}[6cm][l]{\includegraphics[width=6cm]{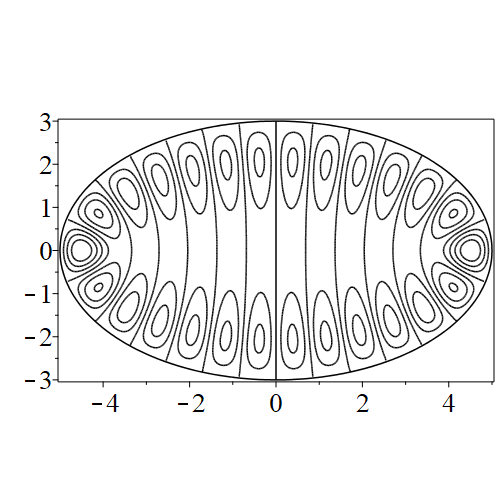}}
  \newline\noindent
  \subcaptionbox{$g=11$\label{fig:cwhisper11}}[6cm][l]{\includegraphics[width=6cm]{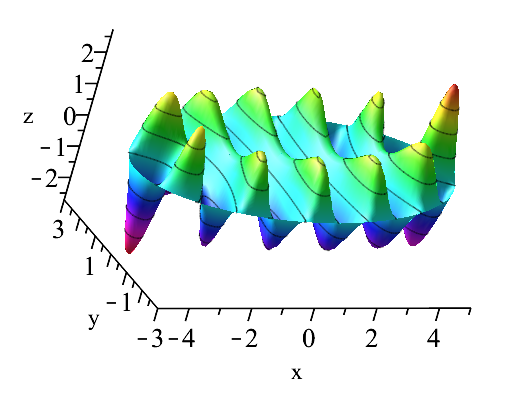}}
  \subcaptionbox{$g=13$\label{fig:cwhisper13}}[6cm][l]{\includegraphics[width=6cm]{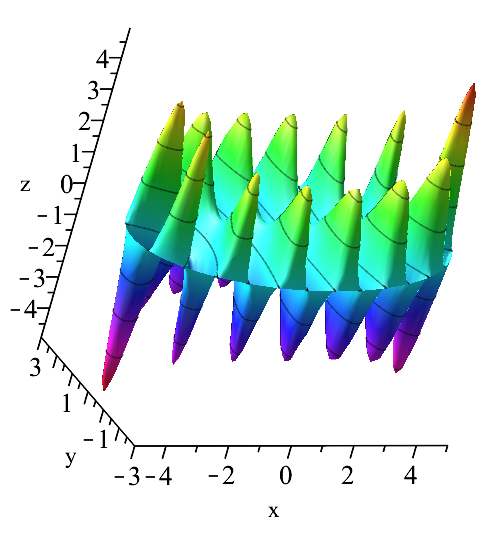}}
  \caption{Two ``whispering gallery'' modes, after~\cite{chen1994visualization}. These are contours of $\Ce_{g}(q;\beta)\ce_{g}(q;\alpha)$ for $g=11$ and $q \approx 56.647$ (left pair) and for $g=13$ with $q \approx 74.437$ (right pair). The name is because most of the vibration occurs near the boundary: the idea is that if these were rooms with elliptic shape, a person at one vertex could whisper at the appropriate frequency and be heard by a person at the other vertex, while another person in the middle could not hear what was being said. The colour plots were produced using \texttt{plots[surfdata]}.\label{fig:whispers}}
\end{figure}

The paper~\cite{chen1994visualization} intriguingly describes ``whispering gallery'' modes, two of which I plot in figure~\ref{fig:whispers}), and ``bouncing ball'' modes, which seem to be exhibited in the (d) figures corresponding to the fourth zero of all modes plotted in figures~\ref{fig:a0}--\ref{fig:b5}.  In a ``bouncing ball'' mode, most of the vibration is confined to a central channel.

\section{Concluding remarks}
None of the pictures in this paper would have surprised Mathieu.  He predicted the hyperbolic and elliptic nodal lines, and the exact numbers of the hyperbolic lines by an ingenious sign variation argument with Sturm sequences.  He also noted that his methods could be used to solve the ``confocal elliptic lake'' problem---that is, a drum with a confocal hole taken out from the middle, although I did not read his description of the necessary boundary conditions all that closely (these boundary conditions are covered briefly in Chapter 28 of the Digital Library of Mathematical Functions~\href{https://dlmf.nist.gov/28}{https://dlmf.nist.gov/28} and more thoroughly in~\cite{gutierrez2003mathieu}).

From the~$q$ data in tables~\ref{tab:qa} and~\ref{tab:qb} we may deduce the frequencies~$\lambda$ of each of these modes.  Writing the general solution as a sum, we have
\begin{align}\label{eq:fullsoln}
  u(x,y,t) =& \sum_{g\ge 0} \sum_{k\ge 1} A_{g,k} \Ce_g(q_{g,k}^e; \beta)\ce_{g}(q_{g,k}^e; \alpha )\sin2\lambda_{g,k}^e m t \nonumber\\
            &{}+ \sum_{g\ge 1} \sum_{k\ge 1} B_{g,k} \Se_g(q_{g,k}^o; \beta)\se_{g}(q_{g,k}^o; \alpha )\sin2\lambda_{g,k}^o m t\>,
\end{align}
where the relation $\lambda_{g,k}^{e,o} = \sqrt{q_{g,k}^{e,o}}/c$ ties each frequency to its mode; we have only tabulated the lowest frequency modes in this paper.  The unknown coefficients $A_{g,k}$ and $B_{g,k}$ must be determined by the initial displacement impulse.

\bibliographystyle{plain}
\bibliography{bib}
\vfill
\newpage

\end{document}